\documentclass{article}

\usepackage{amsfonts}
\usepackage{amssymb}
\usepackage{amsthm}
\usepackage{amsmath}
\usepackage{newlfont}
\usepackage{graphicx}
\usepackage[all,arc]{xy}
\usepackage{epsfig}

\usepackage{multicol}
\usepackage[T1]{fontenc}
\newtheorem {proposition}{Proposition}[section]
\newtheorem {theorem}{Theorem}[section]
\newtheorem {lemma}{Lemma}[section]
\newtheorem {example}{Example}[section]
\newtheorem {definition}{Definition}[section]
\newtheorem {problem}{Problem}

\newtheorem {corollary}{Corollary}[section]

\newtheorem {conjecture}{Conjecture}

\author{{\DJ{}or\dj{}e Barali\' {c}}\\ {\small Mathematical Institute SASA}\\[-2mm] {\small Belgrade, Serbia}
}

\title{On Non-zero Degree Maps between \\Quasitoric
$4$-Manifolds}
\date{}
\begin{document}
\maketitle

\begin{abstract} We study the map degrees between quasitoric
$4$-manifolds. Our results rely on Theorems proved by Duan and
Wang in \cite{Duan1} and \cite{Duan2}. We determine the set $D (M,
N)$ of all possible map degrees from $M$ to $N$ when $M$ and $N$
are certain quasitoric 4-manifolds. The obtained sets of integers
are interesting, e. g. those representable as the sum of two
squares $D (\mathbb{C} P^2\sharp \mathbb{C} P^2, \mathbb{C} P^2)$
or the sum of three squares $D (\mathbb{C} P^2\sharp \mathbb{C}
P^2 \sharp \mathbb{C} P^2, \mathbb{C} P^2)$. Beside the general
results about the map degrees between quasitoric 4-manifolds, the
connections among Duan-Wang's approach, the quadratic forms, the
number theory and the lattices is established.

\end{abstract}

\renewcommand{\thefootnote}{}
\footnotetext{This research was supported by the Grant 174020 of
the Ministry for Education and Science of the Republic of Serbia.}

\section{Introduction}

The \textit{mapping degree} is one of the earliest topological
invariants and almost every textbook has section devoted to the
definition and the calculations of this invariant. For given two
orientable $n$-manifolds $M$ and $N$, every map $f: M\rightarrow
N$ induces the homomorphism $$f_\ast : H_\ast (M)\rightarrow
H_\ast (N).$$ The degree of $f$ is defined as an integer such that
$$f_\ast ([M])=k \,[N],$$ where $[M]\in H_n (M)$ and $[N]\in H_n
(N)$ are the fundamental class of $M$ and $N$ respectively. The
one of fundamental question in manifold topology is:

\begin{problem}\label{problem1} For given two manifolds $M$ and $N$ what are the
integers that could be realized as degree of some map $$f:
M\rightarrow N?$$
\end{problem}

This is natural, but nontrivial and hard problem.

\begin{definition}For given two closed orientable $n$-manifolds $M$ and $N$, $D (M, N)$
is the set of integers that could be realized as degree of a map
from $M$ to $N$ $$D (M, N)=\{\mathrm{deg} f\left|\right.
f:M\rightarrow N\}.$$
\end{definition}

In dimension 2, the problem of determining $D (M, N)$ is
completely solved, \cite{Knes} and \cite{Edmon2}. In dimension 3,
the problem has been studied by several authors and for numerous
classes of $3$-manifolds is solved. The most important results
about $3$-manifolds could be found in the survey article of Wang
\cite{Wang}. However, even in dimension $2$ and $3$ it is seen
that $D (M, N)$ highly depends on the homotopy types of both $M$
and $N$.

From the standard topology course (\cite{Hatcher} and \cite{Bre})
we know several effective methods for calculating  the  mapping
degree. Proposition 2.30 and Exercises 8. p. 258 in \cite{Hatcher}
could be easily generalized and summarized in:

\begin{theorem}\label{funth1} For a map $f:M\rightarrow N$ between connected
closed orientable $n$-manifolds and a point $y\in N$ such that
$f^{-1} (y)=\{x_1, \dots,x_k\}$ and there is ball $B\subset N$,
$y\in B$ such that $f^{-1} (y)$ is the union of $k$ disjoint balls
$B_1$, $\dots$, $B_k$, $x_i\in B_i$ for every $i$, $1\leq i\leq
k$, the mapping degree $\mathrm{deg} f$ is the sum
$$\mathrm{deg} f=\sum_{i=1}^k  \mathrm{deg} f\left|\right. x_i$$
where $\mathrm{deg} f\left|\right. x_i$ is the local map degree,
i. e. the degree of map $f: \partial B_i\rightarrow \partial B$.
\end{theorem}

\begin{figure}[h!h!]
\centerline{\epsfig{figure=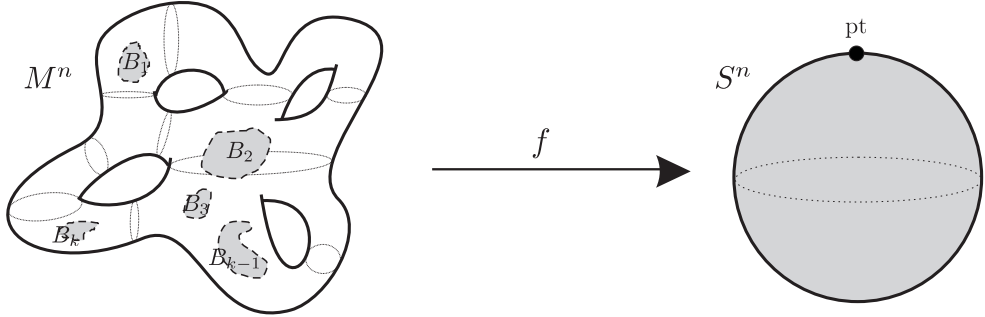,width=\textwidth}}
\caption{The degree $k$ map from $M^n$ to $S^n$.} \label{usferu}
\end{figure}

Theorem \ref{funth1} states that $\mathrm{deg} f$ evaluates the
number of times the domain manifold $M$ "wraps around" the range
manifold $N$ under the mapping $f$. This geometrical principle is
the guiding idea in the most papers studying the mapping degrees.
From Theorem \ref{funth1} it is easy to produce the map of any
given degree into the sphere $S^n$. We take $k$ disjoint balls on
$M$ and map their interiors by orientation preserving
homeomorphism onto $S^n-\{\mathrm{pt}\}$ and the rest of $M$ maps
to the point $\{\mathrm{pt}\}$, Figure \ref{usferu}. Thus, $D
(M^n, S)=\mathbb{Z}$.

Every map $f: M\rightarrow N$ induces homomorphisms on homology
$f_\ast$ and cohomology $f^\ast$. From the following commutative
diagram (see \cite{Hatcher}, p. 241) $$\xymatrix@=20pt{
H_n (M; \mathbb{Z}) \times  \ar@<5ex>[d]^-{f_\ast} H^k (M; \mathbb{Z})  \ar[r]^-{\cap}  & H_{n-k} (M; \mathbb{Z}) \ar[d]^-{f_\ast}\\
H_n (N; \mathbb{Z}) \times H^k (N; \mathbb{Z})
\ar@<5ex>[u]^-{f^\ast} \ar@{>}[r]^-{\cap} & H_{n-k} (M;
\mathbb{Z})}$$ we conclude that for non-zero degree map $f$, every
$f^\ast : H^k (N; \mathbb{Z})\rightarrow H^k (M; \mathbb{Z})$ is
monomorphism. Thus there is no nonzero degree map $f: \Sigma^n
\rightarrow M$ from homology sphere $\Sigma^n$ to the oriented
closed manifold $M$ with nontrivial cohomology.

It is easy to produce the maps of zero degree, so $0\in D (M, N)$.
The identitety map shows that $1\in D (M, M)$. For the beginning,
we reduce Problem \ref{problem1} on:

\begin{problem}\label{problem2} For given two manifolds $M$ and
$N$ is there any degree 1 map from $M$ to $N$?
\end{problem}

The following simple examples show that in general, the sets $D
(M, N)$ are different for distinct manifolds. Even for Problem 2
it is hard to say when the answer is positive or negative.

\begin{proposition} For a simply connected closed orientable manifold
$M^{2n-1}$, the set $D (M, \mathbb{R} P^{2n-1})=2 \mathbb{Z}$.
\end{proposition}

\noindent {\bf Proof:} The sphere $S^n$ is the universal covering
for $\mathbb{R} P^n$ and let $p$ be the covering map. Let $f :
M\rightarrow \mathbb{R} P^n$ be a map. We consider the following
diagram
$$\xymatrix@=20pt{
 & S\ar[d]^-{p}\\
M  \ar@{>}[ur]^-{\tilde{f}} \ar@{>}[r]^-{f} & \mathbb{R} P^{n}.}$$
Because $\pi_1 {M}$ is trivial, there is lifting map $\tilde{f}: M
\rightarrow S^n$ of $f$. From the functoriality of homology we
have $f_\ast = p_\ast {\tilde{f}}_\ast$ and thus $\mathrm{deg}
f=\mathrm{deg} p \cdot \mathrm{deg} \tilde{f}$. But $\mathrm{deg}
\,p=2$ when $n$ is odd and so $\mathrm{deg} f$ is even.

However, it is not hard to produce a degree $2k$ map. Take any
degree $k$ map from $M$ to $S^{2n-1}$ and compose it with $p$.
\hfill $\square$

\begin{proposition}\label{p2} Let $M$ be a simply connected closed orientable $3$-manifold
$M$, and $P$ the Poincar\'{e} homology sphere. Then $D (M, P)=120
\mathbb{Z}$.
\end{proposition}

\noindent {\bf Proof:} Here $S^3$ is the universal cover of $P$
and degree of $p$ is 120. The argument is the same as in the
previous proof. \hfill $\square$

\textit{Remark.} Proposition \ref{p2} is a special case of the
result of Legrand-Matveev-Zieschang \cite{LMZ} about the
calculation of the mapping degree from Seifert manifolds into the
Poincar\'{e} homology sphere. However, this result together with
\cite{LMZ} answered the question for Seifert manifolds.

The results about $3$-manifolds usually suppose some additional
geometrical or topological structure on manifolds. Until recent
years there have been studied some special examples for manifolds
in dimensions above 3, see \cite{Baues}, \cite{Duan3} and
\cite{Gitler}. Haibao Duan and Shicheng Wang progress in the
problem \cite{Duan1} and \cite{Duan2} is significant because it
gives algebraic conditions for the existence of certain map degree
between two given closed $(n-1)$-connected $2n$-manifolds. Their
algebraic conditions are obtained from topology of this wide class
of manifolds. However, even in dimension 4 where the situation is
the simplest, in general it is not easy to check these conditions.

\medskip

The goal of this article is to improve their results in dimension
4 and apply them for the concrete class of $4$-dimensional
quasitoric manifolds. We get interesting results and we put the
results of Duan and Wang in the interaction with recent
development in toric topology and other ideas from topology and
geometry. Quasitoric manifolds are topological generalizations of
toric varieties and represent very important class of manifolds.
They are studied intensively in past twenty years.

In Section $2$ we give review on Duan and Wang work. We prove
again Theorem 2 \cite{Duan1} with accent on some details omitted
in their paper, because we give slightly improved results.

In Section $3$ we briefly discuss topology and geometry of
Hirzebruch surfaces and quasitoric manifolds. In dimension 4 we
have complete classification of quasitoric manifolds. Cohomology
and intersection form of these manifolds is known.

In Section $4$ we study the mapping degrees between
$4$-dimensional quasitoric manifolds.

In Section $5$ we study the mapping degrees between connected sums
of $\mathbb{C} P^2$ and discus the connection of our problem and
the problems about the lattice discriminants.

In Section $6$ we formulate some general theorems about the
mapping degrees between quasitoric 4-manifolds.

In Section $7$ we briefly touch some possible extensions of our
results for higher dimensions and their connections with other
areas of mathematics.

\section{Results of Duan and Wang}

In the articles \cite{Duan1} and \cite{Duan2} are given theorems
which significantly contribute to our knowledge about map degrees
between closed orientable $2n$-manifolds. In this section, we
prove Theorem 2 from \cite{Duan1}, and extend Corollary 3 of
Wang's and Duan's result.

Let $M$ be a $2n$-dimensional closed, connected and orientable
manifold, $n>1$ and let $\bar{H}^n (M; \mathbb{Z})$ be the free
part of $H^n (M; \mathbb{Z})$. Then the cup product operator
$$\bar{H}^n (M; \mathbb{Z})\otimes\bar{H}^n (M; \mathbb{Z})\rightarrow H^{2n} (M; \mathbb{Z})$$
defines the intersection form $X_M$ over $\bar{H}^n (M;
\mathbb{Z})$, which is bilinear and unimodular by Poincar\'{e}
duality, see \cite{Hatcher} Proposition 3.38. This form is
$n$-symmetric in the sense that $$X_M (x\otimes y)=(-1)^n X_M
(y\otimes x).$$

Let $\alpha=(\alpha_1, \dots, \alpha_m)$ be a basis for $\bar{H}^n
(M; \mathbb{Z})$. Then $X_M$ determines an $m\times m$ matrix
$A=(a_{ij})$ where $a_{ij}$ is given by $$a_{i j}=\alpha_i \cup
\alpha_j [M],$$ and $[M]$ is the fundamental class of $H_{2n}
(M)$.

Let $f:M\rightarrow L$ be a map between two connected, closed and
orientable $2n$-manifolds $M$ and $L$, and let $f^\ast$ and
$f_\ast$ be the induced homomorphisms on the cohomology rings and
homology rings. Let $\alpha=(\alpha_1, \dots, \alpha_m)$ and
$\beta=(\beta_1, \dots, \beta_l)$ be basis for $\bar{H}^n (M;
\mathbb{Z})$ and $\bar{H}^n (L; \mathbb{Z})$ respectively. The
induced homomorphism $f^\ast$ determines $m\times l$ matrix
$P=(p_{i j})$ such that $$f^\ast (\alpha_i) = \sum_{j=1}^l p_{i j}
\beta_j,$$ for every $i$, $1\leq i \leq m$.

\begin{theorem}[H. Duan, S. Wang] \label{dwteo}
Suppose $M$ and $L$ are closed oriented $2n$-manifolds with
intersection matrices $A$ and $B$ under some given basis $\alpha$
for $\bar{H}^n (M; \mathbb{Z})$ and $\beta$ for $\bar{H}^n (L;
\mathbb{Z})$. If there is a map $f : M\rightarrow L$ of degree $k$
such that $f^\ast (\beta)=\alpha P$, then $$P^t A P=k B.$$
Moreover, if $k=1$, then $X_L$ is isomorphic to a direct summand
of $X_M$.
\end{theorem}

\noindent {\bf Proof:} For a map $f : M\rightarrow L$ holds
$f_\ast ([M])=k [L]$. From the functoriality of the cup and the
cap product functor we have
$$X_M (f^\ast (x)\otimes f^\ast (y))=f^\ast (x)\cup f^\ast (y)
[M]=$$ $$=f^\ast (x\cup y) [M]=(x\cup y) f_\ast ([M])=(x\cup y)
k[L]=k X_L (x\otimes y)$$ for every $x$, $y\in \bar{H}^n (L;
\mathbb{Z})$. Thus the following diagram commutes
$$\xymatrix@=20pt{
\bar{H}^n (L; \mathbb{Z})\times \bar{H}^n (L; \mathbb{Z})  \ar[r]^-{X_L} \ar[d]^-{f^\ast\otimes f^\ast} & \mathbb{Z} \ar[d]^-{\times k}\\
\bar{H}^n (M; \mathbb{Z}) \times \bar{H}^n (M; \mathbb{Z})
\ar@{>}[r]^-{X_M} & \mathbb{Z}.}$$ Consequently, for the basis
$\alpha$ for $\bar{H}^n (M; \mathbb{Z})$ and $\beta$ for
$\bar{H}^n (L; \mathbb{Z})$ this fact is written in form $P^t A
P=k B$ where $f^\ast (\beta)=\alpha P$.

In particular, when $k=1$ the restriction of $X_M$ on the subgroup
$f^\ast (\bar{H}^n (L; \mathbb{Z}))\subset \bar{H}^n (M;
\mathbb{Z})$ is isomorphic to $X_L$ and unimodular. By Orthogonal
Decomposition Lemma \cite{Milno}, p. 5, $$X_M=X_{f^\ast (\bar{H}^n
(L; \mathbb{Z}))}\oplus X_{H^\bot}=X_L\oplus X_{H^\bot},$$ where
$H^\bot$ is the orthogonal complement of $f^\ast (\bar{H}^n (L;
\mathbb{Z}))$ and $X_{H^\bot}$ is the restriction of $X_M$ on
${H^\bot}$. \hfill $\square$

In the same paper Duan and Wang proved the following theorem that
gives the complete criteria for the existence of degree $k$ map
from a $4$-manifold $M$ to a simply connected 4-manifold $L$.

\begin{theorem}\label{dwt1} Suppose $M$ and $L$ are closed oriented
$4$-manifolds with intersection matrices $A$ and $B$ under given
bases $\alpha$ for $\bar{H}^2 (M; \mathbb{Z})$ and $\beta$ for
$\bar{H}^2 (L; \mathbb{Z})$. If $L$ is simply connected, then
there is a map $f : M \rightarrow L$ of degree $k$ such that
$f^\ast (\beta)= \alpha P$ if and only if $$P^t A P=k B.$$
Moreover there is a map $f:M\rightarrow L$ of degree 1 if and only
if $X_L$ isomorphic to a direct summand of $X_M$.
\end{theorem}

Duan and Wang proved in \cite{Duan1} Corollary 3 of Theorem
\ref{dwteo}. Implicitly, it is clear from their paper this
corollary could be generalized.

\begin{corollary}\label{corb} Suppose $M$ and $L$ are closed oriented
$2n$-manifolds such that \\ $\mathrm{rank} \bar{H}^n (M;
\mathbb{Z})=\mathrm{rank} \bar{H}^n (L; \mathbb{Z})=2 r+1$.
 Then for any map $f: M\rightarrow L$, the absolute value of the degree
of $f$ is a square of an integer.
\end{corollary}

\noindent {\bf Proof:} Let $P$ be the matrix realized by $f$. By
Theorem \ref{dwteo} we have $$P^t A P=k B,$$ where $P$, $A$ and
$B$ are square matrices of order $2 r+1$. By taking the
determinant we have $|P|^2 |A|=k^{2 r+1} |B|$. Since $A$ and $B$
are unimodular, then $|P|^2=|k|^{2 r+1}$. Thus, $|k|$ is a perfect
square. \hfill $\square$

In papers \cite{Duan1} and \cite{Duan2}, Duan and Wang developed
technique for studying the non-zero degree maps between
$(n-1)$-connected closed and oriented $2n$-manifolds. They
demonstrated application on various concrete examples of
manifolds. In Section 4, we continue in the same manner, trying to
put the results in \cite{Duan1} Section 7 in more general context
and classify the obtained results.

\section{Quasitoric manifolds}

Quasitoric manifolds from the general topological point of view
are manifolds with a certain nice torus action. This notion
appeared first in \cite{DaJan}. These objects are connected with
toric varieties from algebraic geometry. In monograph \cite{BuPan}
there is nice exposition about this topic. For purposes of our
results we offer brief overview on them following this exposition.

\subsection{Quasitoric manifolds and characteristic maps}

Let $\mathbb{C}^\ast=\mathbb{C}\setminus \{0\}$ denote the
multiplicative group of complex numbers. The product
$(\mathbb{C}^\ast)^n$ is \textit{algebraic torus}. The torus
$T^n=\underbrace{S^1\times\dots\times S^1}_{n}$ is a subgroup of
the algebraic torus $(\mathbb{C}^\ast)^n$  in the standard way:
$$T^n=\{(z_1, \dots, z_n)\in \mathbb{C}^n \left|\right.
|z_1|=\dots=|z_n|=1\}.$$ The torus $T^n$ acts on $\mathbb{C}^n$ by
standard diagonal multiplication and the quotient of this action
is the positive cone $\mathbb{R}^n_+$.

Let $M^{2n}$ be a $2n$-dimensional manifold with an action of the
torus $T^n$ - \textit{$T^n$-manifold}.

\begin{definition} A standard chart on $M^{2n}$ is a triple $(U,
f, \psi)$ where $U$ is a $T^n$-stable open subset of $M^{2n}$,
$\psi$ is an automorphism of $T^n$, and $f$ is a
$\psi$-equivariant homeomorphism $f:U\rightarrow W$ with some
($T^n$-stable) open subset $W\subset \mathbb{C}^n$ i. e. $f (t
\cdot y)=\psi (t) f (y)$ for all $t\in T^n$ and $y\in U$. Say that
a $T^n$-action on $M^{2n}$ is locally standard if $M^{2n}$ has a
standard atlas, that is, every point of $M^{2n}$ lies in a
standard chart.
\end{definition}

The orbit space for a locally standard action of $T^n$ on $M^{2n}$
is an $n$-dimensional manifold with corner and quasitoric
manifolds correspond to the case when this orbit space is
diffeomorphic, as manifold with corners to a simple polytope
$P^n$.

\begin{definition} Given a combinatorial simple polytope $P^n$, a
$T^n$-manifold $M^{2n}$ is called a quasitoric manifold over $P^n$
if the following two conditions are satisfied:
\begin{enumerate}
    \item the $T^n$ action is locally standard;
    \item there is a projection map $\pi : M^{2n}\rightarrow P^n$
    constant on $T^n$-orbits which maps every $k$-dimensional
    orbit to a point in the interior of a codimension-$k$ face of
    $P^n$, $k=0$, $\dots$, $n$.
\end{enumerate}
\end{definition}

\begin{figure}[h!h!]
\centerline{\epsfig{figure=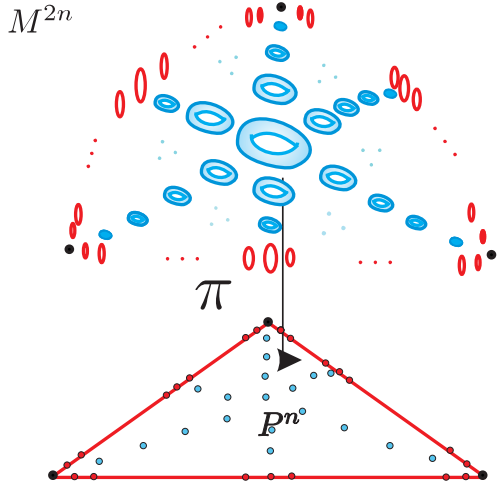,width=0.5\textwidth}}
\caption{The orbit map $\pi$ of quasitoric manifold $M^{2n}$}
\label{quasitoric}
\end{figure}

The $T^n$-action on a quasitoric manifold $M^{2n}$ is free over
the interior of the quotient polytope $P^n$ and vertices of $P^n$
correspond to the $T^n$-fixed points of $M^{2n}$. Let $F_1$,
$\dots$, $F_m$ be facets of $P^n$. For every facet $F_i$, the
pre-image $\pi^{-1} (\mathrm{int} F_i)$ consists of
codimension-one orbits with the same $1$-dimensional isotropy
subgroup, which we denote $T (F_i)$, (see Figure
\ref{quasitoric}). $\pi^{-1} (\mathrm{int} F_i)$ is an $2
(n-1)$-dimensional quasitoric submanifold over $F_i$, with respect
to the action of $T^n/ T(F_i)$ and we denote it $M_i^{2(n-i)}$ and
refer to it as the \textit{facial submanifold} corresponding to
$F_i$. Its isotropy subgroup $T (F_i)$ can be written as $$T
(F_i)=\left\{\left(z^{\lambda_{1i}}, \dots, z^{\lambda_{1
n}}\right)\in T^n\left|\right. |z|=1\right\}.$$ The vector
$\mathbf{\lambda_i}=\left(\lambda_{1i}, \dots, \lambda_{n
i}\right)^t\in  \mathbb{Z}^n$ is determined only up to sign and is
called the \textit{facet vector} corresponding to $F_i$. The
correspondence
\begin{equation}\label{chmap}
l: F_i \mapsto T (F_i)
\end{equation}
is called the \textit{characteristic map} of $M^{2n}$.

Let $G^{n-k}$ be a codimension-$k$ face written as an intersection
of $k$ facets $G^{n-k}=F_{i_1}\cap\dots\cap F_{i_k}$. Then the
submanifolds $M_{i_1}$, $\dots$, $M_{i_k}$ intersect transversally
in a submanifold $M (G)^{2 (n-k)}$, which we refer to as the
\textit{facial submanifold} corresponding to $G$. The map
$T(F_{i_1})\times\dots\times T(F_{i_k})\rightarrow T^n$ is
injective since $T(F_{i_1})\times\dots\times T(F_{i_k})$ is the
$k$-dimensional isotropy subgroup of $M (G)^{2 (n-k)}$. Thus, the
vectors $\mathbf{\lambda}_{i_1}$, $\dots$,
$\mathbf{\lambda}_{i_k}$ form a part of an integral basis of
$\mathbb{Z}^n$.

The correspondence $$G^{n-k}\mapsto \,\mbox{isotropy subgroup
of\,} {M (G)}^{2(n-k)}$$ extends the characteristic map
(\ref{chmap}) to a map from the face poset of $P^n$ to the poset
of subtori of $T^n$.

\begin{definition} Let $P^n$ be a combinatorial simple polytope
and $l$ is a map from facets of $P^n$ to one-dimensional subgroups
of $T^n$. Then $(P^n, l)$ is called a characteristic pair if $l
(F_{i_1})\times\cdots\times l (F_{i_k})\rightarrow T^n$ is
injective whenever $F_{i_1}\cap\dots\cap F_{i_k}\neq\emptyset$.
\end{definition}

The map $l$ directly extends to a map from the face poset of $P^n$
to the poset of subtori of $T^n$, so we have subgroup $l
(G)\subset T^n$ for every face $G$ of $P^n$. As in the case of
standard action of $T^n$ on $\mathbb{C}^n$, there is projection
$T^n\times P^n\rightarrow M^{2n}$ whose fibre over $x\in M^{2n}$
is the isotropy subgroup of $x$. This argument we use for
reconstructing the quasitoric manifold from any given
characteristic pair $(P^n, l)$.

\medskip

Given a point $q\in P^n$, we denote by $G (q)$  the minimal face
containing $q$ in its relative relative. Define relation $\sim$ on
$T^n \times P^n$ in following way $(t_1, p)\sim (t_2, q)$ if and
only if $p=q$ and $t_1 t_2^{-1}\in l (G (q))$. Now set $$M^{2n}
(l):=(T^n\times P^n)/\sim .$$ The free action of $T^n$ on
$T^n\times P^n$ obviously descends to an action on $(T^n\times
P^n)/\sim$, with quotient $P^n$. The latter action is free over
the interior of $P^n$ and its fixed points are vertices of $P^n$,
(see Figure \ref{quasitoric}). Just as $P^n$ it is covered by the
open sets $U_v$, based on the vertices and diffeomorphic to
$\mathbb{R}^n_+$, so the space $(T^n\times P^n)/\sim$ is covered
by open sets $(T^n \times U_v)/\sim$ homeomorphic to $(T^n\times
\mathbb{R}^n_+)/\sim$, and therefore to $\mathbb{C}^n$. This
implies that the $T^n$-action on $(T^n\times P^n)/\sim$ is locally
standard, and therefore $(T^n\times P^n)/\sim$ is a quasitoric
manifold.

\subsection{Cohomology of quasitoric manifolds}

The cohomology of quasitoric manifolds is determined by Davis and
Januszkiewicz \cite{DaJan}. They constructed also $C W$ structure
on quasitoric manifolds with only even dimensional cells.

Let $F_1$, $\dots$, $F_m$ be facets of a simple polytope $P^n$ and
let $\mathbb{Z} [v_1, \dots, v_m]$ be the polynomial algebra on
$\mathbb{Z}$ with $m$ generators $v_1$, $\dots$, $v_m$ one for
each facet. The \textit{Stanley-Reisner ring} of a simple polytope
$P^n$ is the quotient ring  $\mathbb{Z} (P^n)=\mathbb{Z} [v_1,
\dots, v_m]/ \mathcal{I}_P$ where $\mathcal{I}_P$ is the ideal
generated by all square-free monomials $v_{i_1} v_{i_2}\cdots
v_{i_s}$ such that $F_{i_1}\cap\cdots\cap F_{i_s}=\emptyset$ in
$P$, $i_1<\cdots< i_s$.

Given a quasitoric manifold $M^{2n}$ with characteristic map $l:
F_i\mapsto T (F_i)$ and facet vectors
$\mathbf{\lambda_i}=(\lambda_{1i}, \dots, \lambda_{n i})^t\in
\mathbb{Z}^n$, $i=1$, $\dots$, $m$, define linear forms
\begin{equation}\label{ideal:l}
\theta_i:=\lambda_{i 1} v_1+ \cdots, \lambda_{i m} v_m\in
\mathbb{Z} [v_1, \dots, v_m], \quad 1\leq i\leq n.
\end{equation}
The images of these linear forms in the Stanley-Reisner ring
$\mathbb{Z} (P^n)$ will be denoted by the same letters. Let
$\mathcal{I}_l$ denote the ideal in $\mathbb{Z} (P^n)$ generated
by $\theta_1$, $\dots$, $\theta_n$.

\begin{theorem}[Davis and Januszkiewicz] The cohomology ring of
$M^{2n}$ is given by $$H^\ast (M^{2n}; \mathbb{Z})=\mathbb{Z}
[v_1, \dots, v_m]/(\mathcal{I}_P+\mathcal{I}_l)=\mathbb{Z}
(P^n)/\mathcal{I}_l,$$ where $v_i$ is the $2$-dimensional
cohomology class dual to the facial submanifold $M_i^{2 (n-1)}$,
$i=1$, $\dots$, $m$.
\end{theorem}

\subsection{Hirzebruch surfaces and $4$-dimensional quasitoric
manifolds}

\textit{Hirzebruch surfaces} were introduces by Hirzebruch in
\cite{Hirz} and they are algebraic surfaces over the complex
numbers. As complex manifolds they are pairwise distinct while as
smooth manifolds there are only two diffeomorphism types.

Given  an integer $k$, the Hirzebruch surface is the complex
manifold $\mathbb{C} P (\zeta_k \oplus \mathbb{C})$, where
$\zeta_k$ is the complex line bundle over $\mathbb{C} P^1$ with
first Chern class $k$, and $\mathbb{C} P (\cdot)$ denotes the
projectivisation of a complex bundle. Each Hirzebruch surface is
the total space of the bundle $H_k \rightarrow \mathbb{C} P^1$
with fibre $\mathbb{C} P^1$. For even $k$ the surface $H_k$ is
diffeomorphic to $S^2\times S^2$ and for odd $k$ to $\mathbb{C}
P^2\sharp \overline{\mathbb{C} P^2}$, where $\overline{\mathbb{C}
P^2}$ denotes the space $\mathbb{C} P^2$ with reversed
orientation. The Hirzebruch surface $H_k$ is quasitoric manifold
which orbit space is a combinatorial square.

The topological classification problem (up to diffeomorphism) for
quasitoric manifolds over a given simple polytope is intractable.
Up to now, only some particular results are known. In \cite{Orli}
classification problem for quasitoric manifolds over polygons is
completely solved.

\begin{theorem}\label{klasi} A quasitoric manifold of dimension $4$
is diffeomorphic to connected sum of several copies of $\mathbb{C}
P^2$, $\overline{\mathbb{C} P^2}$ and $S^2\times S^2$.
\end{theorem}

Quasitoric manifolds are simply connected, so they are perfect
test examples for the application of Theorem \ref{dwt1}.

\section{Mapping degrees between Quasitoric $4$-manifolds}

The mapping degree between some quasitoric $4$-manifolds are
studied in \cite{Duan1}, where some simply examples illustrated
application of Theorems \ref{dwteo} and \ref{dwt1}. However, the
results they got in simplest case are not trivial. In this section
we focus our attention on more general cases of quasitoric
$4$-manifolds.

Quasitoric $4$-manifolds are topologically classified by Theorem
\ref{klasi} and we could easily determine their intersection form.
The matrices representing the intersection form for  $\mathbb{C}
P^2$, $\overline{\mathbb{C} P^2}$ and $S^2\times S^2$ are
\begin{equation}\label{osnovno}
\left[\begin{array}{c} 1\end{array}\right], \qquad
\left[\begin{array}{c}
-1\end{array}\right],\qquad \mbox{abd} \qquad \left[\begin{array}{cc} 0 & 1\\
1 & 0\end{array}\right]
\end{equation}
respectively. Thus, the intersection form for $\mathbb{C}
P^2\sharp \overline{\mathbb{C} P^2}$ has matrix representation
\begin{equation}\label{osnovno1}\left[\begin{array}{cc} 1 & 0\\
0 & -1\end{array}\right].
\end{equation}

Theorem \ref{klasi} states that a quasitoric $4$-manifold is
diffeomorphic to $$\left(\mathbb{C} P^2\right)^{\sharp a}\sharp
\left(\overline{\mathbb{C} P^2}\right)^{\sharp b}\sharp
\left(S^2\times S^2\right)^{\sharp c}=$$ $$
 =\underbrace{\mathbb{C} P^2\sharp\cdots\sharp\mathbb{C}
P^2}_{a \, \mbox{times}} \sharp \underbrace{\overline{\mathbb{C}
P^2} \sharp\cdots \sharp\overline{\mathbb{C} P^2}}_{b \,
\mbox{times}}\sharp \underbrace{S^2\times S^2 \sharp \cdots \sharp
S^2\times S^2}_{c \, \mbox{times}}.$$ The intersection form of a
given quasitoric manifold has representation in $(a+b+2 c)$-square
matrix
\begin{equation}\label{forma1}
\left[\begin{array}{ccc}
I_{a \times a} & 0 & 0\\
0 & - I_{b\times b} & 0\\
0 & 0 &  A_{c\times c}
\end{array}\right],
\end{equation} where $A_{c \times c}$ is $2c$-square matrix $$\left[\begin{array}{ccccccc} 0 & 1 & 0 & 0 & \cdots & 0 & 0\\
1 & 0 & 0 & 0 & \cdots & 0 & 0\\
0 & 0 & 0 & 1 & \cdots & 0 & 0\\
0 & 0 & 1 & 0 & \cdots & 0 & 0\\
\vdots & \vdots & \vdots & \vdots & \ddots & \vdots & \vdots\\
0 & 0 & 0 & 0 & \cdots & 0 & 1\\
0 & 0 & 0 & 0 & \cdots & 1 & 0\end{array}\right].$$

Let $M=\left(\mathbb{C} P^2\right)^{\sharp a}\sharp
\left(\overline{\mathbb{C} P^2}\right)^{\sharp b}\sharp
\left(S^2\times S^2\right)^{\sharp c}$ and $N=\left(\mathbb{C}
P^2\right)^{\sharp d}\sharp \left(\overline{\mathbb{C}
P^2}\right)^{\sharp e}\sharp \left(S^2\times S^2\right)^{\sharp
f}$ be two quasitoric manifolds, and the matrices $A$ and $B$
their intersection forms respectively. From Theorem \ref{dwt1}
follows that there is degree $k$ map between $M$ and $N$ if and
only if there is $(a+b+2c)\times (d+e+2 f)$ matrix $P$ such that
$P^t A P=k B$. Direct calculation gives that elements of
$(d+e+2f)$-square matrix $C=P^t A P$ are $$c_{i j}=\sum_{r=1}^a
p_{r i} p_{r j}- \sum_{r=1}^b p_{a+r\, i} p_{a+r\, j}+\sum_{r=1}^c
(p_{a+b+2r-1\, i} p_{a+b+2 r\, j}+ p_{a+b+2r\, i} p_{a+b+2r-1\,
j}).$$ Solving the equation $P^t A P= k B$ is then equivalent to
solving the certain system of Diophantine equations. There is no
algorithm for solving a Diophantine equation, so there is not
natural way to approach this problem. However, only $d+e+f$ out of
$\frac{(d+e+f)(d+e+f+1)}{2}$ expressions is equal to $\pm k$ and
the others are $0$.

We prove the following lemma:

\begin{lemma}\label{simpa} Let $B_{n\times k}$ and $C_{k\times n}$ be arbitrary
matrices such that $k<n$. Then $$\det |B\cdot C|=0.$$
\end{lemma}

\noindent {\bf Proof:} Denote $A=B\cdot C$. Then $a_{i
j}=\sum_{r=1}^k b_{i r} c_{r k}$. Let $B'$ and $C'$ be $n$ square
matrices such that $$B'=\left[\begin{array}{cc} B &
0\end{array}\right],  \qquad C'=\left[\begin{array}{c} C \\
0\end{array}\right]. $$ It is obviously that $$B'\cdot C'=B\cdot
C=A,$$ and consequently $$0=\det B' \cdot \det C'=\det A.$$ \hfill
$\square$

\begin{theorem}\label{b1} Let $M$ and $N$ be closed, connected and oriented
$2n$-manifolds such that $1\leq \mathrm{rank} \bar{H}^n (M;
\mathbb{Z})< \mathrm{rank} \bar{H}^n (N; \mathbb{Z})$. Then there
is no non-zero degree map $f: M\rightarrow N$.
\end{theorem}

\noindent {\bf Proof:} This is corollary of Lemma \ref{simpa}
applied for matrices $P^t$ and $A \cdot P$. According to Lemma
$$0=k^{\mathrm{rank} \bar{H}^n (N; \mathbb{Z})} \det B$$ and
because of unimodularity of the intersection form $k=0$. \hfill
$\square$

 \begin{theorem}\label{pros1} If there is a degree $k$ map $f: M\rightarrow N$
 between two $4$-manifolds, then there is a degree $k$ map from
 $M\sharp \mathbb{C} P^2$ to $N$.
 \end{theorem}

\noindent {\bf Proof:} Let $A$ be the intersection matrix for $M$,
$B$ for $N$ and $P$ matrix such that $$P^t A P= k B.$$ Then
$A'=\left[\begin{array}{cc} 1 & 0\\
0 & A
\end{array}\right]$ is the intersection matrix for  $M\sharp \mathbb{C}
P^2$. The matrix $P'=\left[\begin{array}{c} 0\\
P
\end{array}\right]$ clearly satisfies $$(P')^t A' P'=k B.$$\hfill
$\square$

In the same fashion we can prove the following theorem:

\begin{theorem}\label{pros2} If there is a degree $k$ map $f: M\rightarrow N$
 between two $4$-manifolds, then there is a degree $k$ map from
 $M\sharp \left(S^2\times S^2\right)$ to $N$.
 \end{theorem}

\begin{corollary}\label{bara2} If there is a degree $k$ map $f: M\rightarrow N$
between two $4$-manifolds, then there is a degree $k$ map from
$M\sharp Q$ to $N$ for every $4$-manifold $Q$.
\end{corollary}

\begin{theorem}\label{bar1}  If there are degree $k$ maps $f: M\rightarrow N$
and $g: M'\rightarrow N'$ between  $4$-manifolds, then there is a
degree $k$ map from $M\sharp M'$ to $N\sharp N'$.
 \end{theorem}

\noindent {\bf Proof:} Let $A$ and $A'$ be the intersection
matrices for $M$ and $M'$, $B$ and $B'$ for $N$ and $N'$, $P$ and
$P'$ matrices such that
$$P^t A P= k B,$$  $$P'^t A' P'= k B'.$$ We easily check that
$$\left[\begin{array}{cc}
P^t & 0\\
0 & P'^t\end{array}\right] \cdot \left[\begin{array}{cc}
A & 0\\
0 & A'\end{array}\right] \cdot \left[\begin{array}{cc}
P & 0\\
0 & P'\end{array}\right]=k \left[\begin{array}{cc}
B & 0\\
0 & B'\end{array}\right].$$\hfill $\square$

 \subsection{Maps to $\mathbb{C} P^2$}

 We study maps from quasitoric manifolds to $\mathbb{C} P^2$  (and $\overline{\mathbb{C}
 P^2}$). Let $M$ be a quasitoric manifold diffeomorphic to $\left(\mathbb{C} P^2\right)^{\sharp a}\sharp
\left(\overline{\mathbb{C} P^2}\right)^{\sharp b}\sharp
\left(S^2\times S^2\right)^{\sharp c}$.

Theorem \ref{dwt1} reduces problem on the existence of nontrivial
solution of Diophantine equation \begin{equation}\label{CP2}
\sum_{i=1}^a p^2_{i \,1}- \sum_{i=1}^b p_{a+i \, 1}^2+2
\sum_{i=1}^c p_{a+b+2 i-1\, 1} p_{a+b+2 i\, 1}=k.\end{equation}

\begin{theorem}\begin{itemize}
    \item If $a\geq 1$ (or $b\geq 1$) and $c\geq 1$  then the equation \ref{CP2}
    have solution for every $k\in \mathbb{Z}$
    \item If $a\geq 4$ and $b=c=0$ then the equation \ref{CP2}
    have solution for every nonnegative integer $k$ and no
    solution for negative $k$.
    \item If $b\geq 4$ and $a=c=0$ then the equation \ref{CP2}
    have solution for every integer $k\leq 0$ and no
    solution for positive $k$.
    \item If $a=3$ and $b=c=0$ then the equation \ref{CP2}
    have solution for every nonnegative integer $k\neq 4^p (8
    q+7)$ and no solution for positive integers $k=4^p (8
    q+7)$ and negative integers.
    \item If $b=3$ and $a=c=0$ then the equation \ref{CP2}
    have solution for every integer $k\neq -4^p (8
    q+7)$ and no solution for negative integers $k=-4^p (8
    q+7)$ and positive integers.
    \item If $a=2$ and $b=c=0$ then the equation \ref{CP2}
    have solution for every nonegative integer $k$ such that every
    prime number $4p-1$ that divides $k$ occurs even time in the prime
    factorization of $k$ and no solutions in other cases
    \item If $b=2$ and $a=c=0$ then the equation \ref{CP2}
    have solution for every integer $k\leq 0$ such that every
    prime number $4p-1$ that divides $|k|$ occurs even time in the prime
    factorization of $|k|$ no solutions in other cases.
    \item If $a=b=1$ and $c=0$ then the equation \ref{CP2}
    have solution for every integer $k\neq 4 p+2$ and no solution
    for $k=4 p+2$
    \item If $a=1$ and $b=c=0$ then \ref{CP2}
    have solution for every integer that is square of an integer
    and no solution in other cases
    \item If $b=1$ and $a=c=0$ then \ref{CP2}
    have solution for every integer that is square of an integer
    multiplied by $-1$ and no solution in other cases
    \item If $a=b=0$ and $c\geq 1$ then the equation \ref{CP2}
    have solution for every even integer $k$
\end{itemize}
\end{theorem}

\noindent {\bf Proof:} We observe that every integer could be
written as $u^2+2 v w $ where $u$, $v$ and $w$ are integers. This
guaranties existence of $k$ degree map $f: \mathbb{C} P^2\sharp
(S^2\times S^2)\rightarrow \mathbb{C} P^2$. Theorems \ref{pros1}
and \ref{pros2} extend the result for cases $a\geq 1$ and $c\geq
1$. Since every nonnegative integer has decomposition into the sum
of four squares of integers, there is $k\geq 0$ degree map when
$a\geq 4$ and $b=c=0$.

The case $a=3$ and $b=c=0$ is interesting, it follows from the
nontrivial result of Legendre (1798) and Gauss (1801) that
nonnegative integer has presentation as the sum of three squares
of integers iff it is not of the type $4^p (8 q+7)$.

If $a=2$ and $b=c=0$ then for any map $f: \mathbb{C} P^2\sharp
\mathbb{C} P^2\rightarrow \mathbb{C} P^2$, $\mathrm{deg} f$ must
be the sum of two squares. It is well known fact that nonnegative
integer $k$ could be written as the  sum of two squares iff that
every prime number $4p-1$ that divides $k$ occurs even time in the
prime factorization of $k$. For $k=u^2+v^2$ we take the matrix
$$P=\left[\begin{array}{cc}
u & v\\
v & -u\end{array}\right]$$ and get the map of degree $k$.

In the case $a=b=1$ and $c=0$, the map degree must be the
difference of two squares and this is possible iff the integer is
not equal 2 modulo 4. Again, it is easy to check
$$P=\left[\begin{array}{cc}
u & v\\
v & -u\end{array}\right]$$ gives map of degree $u^2-v^2$.

In the same manner could be proved other cases so we omit the rest
of proof. \hfill $\square$

\subsection{Maps to $S^2\times S^2$}

Maps from $\mathbb{C}P^2\sharp \mathbb{C}P^2$ and
$\mathbb{C}P^2\sharp \overline{\mathbb{C}P^2}$ to $S^2\times S^2$
are studied in \cite{Duan1}. We are interested in the mapping
degrees from an arbitrary quasitoric manifold to $S^2 \times S^2$.

\begin{proposition} There is no non-zero degree map from $({\mathbb{C}P^2})^{\sharp k}$ ( $({\overline{\mathbb{C}P^2}})^{\sharp k}$) to $S^2\times
S^2$.
\end{proposition}

\begin{proposition} For every integer $k$ there is a degree $k$ map from $({S^2\times S^2})^{\sharp k}$ to $S^2\times S^2$.
\end{proposition}

\begin{proposition} For every integer $k$ there is a degree $k$ map from $({\mathbb{C}P^2})^{\sharp 2}\sharp \overline{{\mathbb{C}P^2}} $ \\
$({\mathbb{C}P^2}\sharp (\overline{{\mathbb{C}P^2}})^{\sharp 2})$
to $S^2\times S^2$.
\end{proposition}

\begin{theorem} Let $M$ be a quasitoric manifold such that $c\geq
1$, then there is a degree $k$ map from $M$ to $S^2 \times S^2$
for every integer $k$.
\end{theorem}

\noindent {\bf Proof:} We easily check that matrix $P$ such that
$$P=\left[\begin{array}{c}
\begin{array}{c} 0
\end{array}\\
\begin{array}{cc} 0 & 1\\
k & 0
\end{array}
\end{array}\right],$$ satisfy Theorem \ref{dwt1}. \hfill $\square$

\begin{theorem} Let $M$ be a quasitoric manifold such that $a\geq
2$ and $b\geq 1$ (or $a\geq 1$ and $b\geq 2$) then then there is a
degree $k$ map from $M$ to $S^2 \times S^2$ for every integer $k$.
\end{theorem}

\noindent {\bf Proof:} We easily check that matrix $P$ such that
$$P=\left[\begin{array}{c}
\begin{array}{c} 0_{(a-2)\times 2}
\end{array}\\
\begin{array}{cc} k & 0\\
0 & 1\\
-k & 1
\end{array}\\
0_{(b+2 c-1) \times 2}
\end{array}\right],$$ satisfy Theorem \ref{dwt1}. \hfill $\square$

\subsection{Maps to $\mathbb{C} P^2\sharp \mathbb{C} P^2$}

According to \cite{Duan1} there is no nonzero degree map from
$\mathbb{C} P^2\sharp \overline{\mathbb{C} P^2}$ and $S^2\times
S^2$ to $\mathbb{C} P^2\sharp \mathbb{C} P^2$. For other
quasitoric manifolds the sets $D (M, \mathbb{C} P^2\sharp
\mathbb{C} P^2)$ are richer and interesting.

\begin{proposition}\label{cp22}There is a
degree $k$ map $f: \left(\mathbb{C} P^2\right)^{\sharp 2}
\rightarrow \left(\mathbb{C} P^2\right)^{\sharp 2}$ if and only if
$k\geq 0$ and every prime number $4p-1$ that divides $k$ occurs
even time in the prime factorization of $k$.
\end{proposition}

\begin{proposition}There is a
degree $k$ map $f: (S^2\times S^2)^{\sharp n} \rightarrow
\left(\mathbb{C} P^2\right)^{\sharp 2}$ $n\geq 2$ if and only if
$k$ is even number.
\end{proposition}

\begin{proposition}For every integer $k$ there is a
degree $k$ map \\ $f: \mathbb{C} P^2\sharp \overline{\mathbb{C}
P^2}\sharp(S^2\times S^2)\rightarrow \left(\mathbb{C}
P^2\right)^{\sharp 2}$.
\end{proposition}

\begin{theorem} There is a
degree $k$ map $f: \left(\mathbb{C} P^2\right) \sharp (S^2\times
S^2)\rightarrow \mathbb{C} P^2\sharp \mathbb{C} P^2$ if and only
if $k$ is the square of an integer or the twice square of an
integer.
\end{theorem}

\begin{theorem} There is $k$ degree map $f: \left(\mathbb{C} P^2\right)^{\sharp 2} \sharp
(S^2\times S^2)\rightarrow \mathbb{C} P^2\sharp \mathbb{C} P^2$ if
and only if $k\geq 0$ and every prime number $4p-1$ that divides
$k$ occurs even time in the prime factorization of $k$.
\end{theorem}

\noindent {\bf Proof:} Let $$P=\left[
\begin{array}{cc}a & b\\
c & d\\
e & f\\
g & h
\end{array}\right].$$ We are solving the system \begin{eqnarray}
a^2+c^2+2 g e =b^2+d^2+2 f h=k\\ ab+cd+ef+g h=0.
\end{eqnarray}
It is clear that when $k=m^2+n^2$ we have solution $a=d=m$, $b=n$,
$c=-n$ and $g=e=h=k=0$.

We are going to prove that $k$ has the form $m^2+n^2$.

$$(k-(a^2+c^2))(k-(b^2+d^2))=4 g e f h= 2(a b+c d)^2-2 (g e)^2-(f
h)^2$$
$$k^2+(bc-a d)^2+(g h-e f)^2=k (a^2+b^2+c^2+d^2).$$

It is clear that $k\geq 0$. Let $k$ be minimal integer such there
is prime number $q=4 p-1$ such that $(q^{2 r+1}\mid\mid k$ and
there exist the solution of system. Then $q^{2 r+1}\mid (b c-a
d)^2+(g h-e f)^2$ and consequently $q^{r+1}\mid bc-ad$ and
$q^{r+1}\mid g h-e f$. Thus, $q^{2 r+2}$ divides the left side,
and we have that $q\mid a^2+b^2+c^2+d^2$. Combining this with
$q\mid ad-bc$ we get $q| (a-d)^2+(b+c)^2$ and $q|
(a+d)^2-(b-c)^2$. Thus, $a\equiv b \equiv c\equiv d\equiv 0
\,\mathrm{mod} \,q$ and $q| g e$, $q|g h$. So there exist the
solution for $\frac{k}{q^2}$, contradicting the minimality of $k$.
This proves $k$ could be written as the sum of two squares of
integers. \hfill $\square$

Using the same approach we obtain the following theorems. We omit
the proofs because they are more or less analogous to the proofs
of the previous statements.

\begin{theorem} There is a
degree $k$ map $f: \left(\mathbb{C} P^2\right)^{\sharp 3}
\rightarrow \mathbb{C} P^2\sharp \mathbb{C} P^2$ if and only if
$k\geq 0$ and every prime number $4p-1$ that divides $k$ occurs
even time in the prime factorization of $k$.
\end{theorem}

\begin{theorem} There is a
degree $k$ map $f: \left(\mathbb{C} P^2\right)^{\sharp n}
\rightarrow \mathbb{C} P^2\sharp \mathbb{C} P^2$ $n\geq4$, if and
only if $k$ is nonnegative integer.
\end{theorem}

\begin{theorem} There is a
degree $k$ map $f: \left(\mathbb{C} P^2\right)^{\sharp 3} \sharp
(S^2\times S^2)\rightarrow \mathbb{C} P^2\sharp \mathbb{C} P^2$
for every integer $k$.
\end{theorem}

\begin{theorem} There is a
degree $k$ map $f: \left(\mathbb{C} P^2\right)^{\sharp
2}\sharp\overline{\mathbb{C} P^2} \rightarrow \mathbb{C} P^2\sharp
\mathbb{C} P^2$ if and only if $k\geq 0$ and every prime number
$4p-1$ that divides $k$ occurs even time in the prime
factorization of $k$.
\end{theorem}

\begin{theorem} There is a
degree $k$ map $f: \left(\mathbb{C} P^2\right)^{\sharp
3}\sharp\overline{\mathbb{C} P^2} \rightarrow \mathbb{C} P^2\sharp
\mathbb{C} P^2$ if and only if $k\geq 0$ and every prime number
$4p-1$ that divides $k$ occurs even time in the prime
factorization of $k$.
\end{theorem}

\begin{theorem} There is a
degree $k$ map $f: \left(\mathbb{C} P^2\right)^{\sharp
n}\sharp\overline{\mathbb{C} P^2} \rightarrow \mathbb{C} P^2\sharp
\mathbb{C} P^2$, $n\geq 4$ if and only if $k\geq 0$.
\end{theorem}

\begin{theorem} There is $k$ degree map $f: \left(\mathbb{C} P^2\right)^{\sharp 2}\sharp(\overline{\mathbb{C} P^2})^{\sharp 2} \rightarrow \mathbb{C} P^2\sharp \mathbb{C} P^2$ for every integer $k$.
\end{theorem}

\begin{corollary} If $M$ is a quasitoric $4$-manifold such that
$\mathrm{rank} \bar{H}^2 (M; \mathbb{Z})\geq 5$ and $b+2 c\geq 2$
then for every integer $k$ there is a degree $k$ map
$f:M\rightarrow \mathbb{C} P^2\sharp \mathbb{C} P^2$.
\end{corollary}
\subsection{Maps to $\mathbb{C} P^2\sharp \overline{\mathbb{C} P^2}$}

\begin{proposition} There is no non-zero degree map $f:{\mathbb{C} P^2}^{\sharp n}\rightarrow\mathbb{C} P^2\sharp \overline{\mathbb{C}
P^2}$.
\end{proposition}

\begin{proposition} There is a
degree $k$ map $f:{\mathbb{C} P^2}\sharp \overline{\mathbb{C}
P^2}\rightarrow\mathbb{C} P^2\sharp \overline{\mathbb{C} P^2}$ iff
$k\neq 4 t+2$.
\end{proposition}

\begin{theorem} For every integer $k$ there is a
degree $k$ map $f: ({\mathbb{C} P^2})^{\sharp 2}\sharp
\overline{\mathbb{C} P^2}\rightarrow\mathbb{C} P^2\sharp
\overline{\mathbb{C} P^2}$.
\end{theorem}

\noindent {\bf Proof:} Let $$P=\left[
\begin{array}{cc}a & b\\
c & d\\
e & f
\end{array}\right].$$ We are looking for the solutions of \begin{eqnarray}
a^2+c^2-e^2 =f^2-b^2-d^2=k\\ ab+cd-ef=0.
\end{eqnarray}

For $k\neq 4 t+2$ it is known that there are integers $m$ and $n$
such that $k=m^2-n^2$. In this case, $a=b=0$, $c=f=m$ and $d=e=n$
finishes job. For $k=4 t+2$, we could take $a=1$, $b=2$, $c=2
t+1$, $d=2 t+2$, $e=2 t$ and $f=2 t+3$. \hfill $\square$

\begin{theorem} For every integer $k$ there is a
degree $k$ map $f:{\mathbb{C} P^2}\sharp (S^2\times
S^2)\rightarrow\mathbb{C} P^2\sharp \overline{\mathbb{C} P^2}$.
\end{theorem}

\noindent {\bf Proof:} Let $$P=\left[
\begin{array}{cc}a & b\\
c & d\\
e & f
\end{array}\right].$$ We are looking for the solutions of the system \begin{eqnarray}
a^2+2 c e =b^2+2 d f=k\\ ab+cf+de=0.
\end{eqnarray}

For $k= 2t$ we could take $a=b=0$, $c=d=t$, $e =1$,and $f=-1$. For
$k=2 t+1$, we could take $a=b=d=e=2 t+1$, $e=- t$ and $f=-t-1$.
\hfill $\square$

\begin{corollary} For every quasitoric 4-manifold $M$ such that
$b\geq 1$ or $c\geq 1$ and every integer $k$ there is a degree $k$
map $f: M\rightarrow\mathbb{C} P^2\sharp \overline{\mathbb{C}
P^2}$.
\end{corollary}

\section{Orthogonal lattices and maps between connected sums of $\mathbb{C} P^2$}

In this section we are focused on the maps between connected sums
of $\mathbb{C} P^2$. Our main interest is in the mapping degrees
$$f:({\mathbb{C} P^2})^{\sharp \, n} \rightarrow ({\mathbb{C} P^2})^{\sharp \,
n}.$$

\begin{proposition} There is a degree $k$ map $f:({\mathbb{C} P^2})^{\sharp \, 2n-1} \rightarrow ({\mathbb{C} P^2})^{\sharp \,
2n-1}$, $n\geq 1$ if and only if $k$ is the square of an integer.
\end{proposition}

\begin{theorem}\label{cp4} There is a degree $k$ map $f:({\mathbb{C} P^2})^{\sharp \, 4} \rightarrow ({\mathbb{C} P^2})^{\sharp \,
4}$  if and only if $k$ is nonnegative integer.
\end{theorem}

\noindent {\bf Proof:} We use the fact that every nonnegative
integer could be written as the sum of four squares
$$k=a^2+b^2+c^2+d^2.$$ Then the matrix $$P=\left[
\begin{array}{cccc}a & b & c & d\\
b & -a & -d & c\\
c & d & -a & -b\\
d & -c & b & -a
\end{array}\right]$$ guaranties the existence of a degree $k$ map.\hfill $\square$

Theorem \ref{cp4} together with Theorem \ref{bar1} implies that:

\begin{corollary} There is a degree $k$ map $f:({\mathbb{C} P^2})^{\sharp \, 4 n} \rightarrow ({\mathbb{C} P^2})^{\sharp \,
4 n}$, $n\geq 1$ if and only if $k$ is nonnegative integer.
\end{corollary}

The remaining case to determine the mapping degrees
$f:({\mathbb{C} P^2})^{\sharp \, 4 n+2} \rightarrow ({\mathbb{C}
P^2})^{\sharp \, 4 n+2}$, $n\geq 1$ is still open. Proposition
\ref{cp22} implies that the set of all integers that could be
written as the sum of two squares belongs to $D (({\mathbb{C}
P^2})^{\sharp \, 4 n+2}, ({\mathbb{C} P^2})^{\sharp \, 4 n+2})$.
We could not give the answer even in the case $f:({\mathbb{C}
P^2})^{\sharp \, 6} \rightarrow ({\mathbb{C} P^2})^{\sharp \, 6}$,
but we checked directly that there is no degree $3$, $7$, $11$,
$15$, $19$ and various other cases that afford conjecture that $D
(({\mathbb{C} P^2})^{\sharp \, 6}, ({\mathbb{C} P^2})^{\sharp \,
6})$ is the set of integers that could be written as the sum of
squares. Generally, we suppose:

\begin{conjecture}\label{conj1} The set $D (({\mathbb{C} P^2})^{\sharp \, 4 n+2},
({\mathbb{C} P^2})^{\sharp \, 4 n+2})$ is the set of nonegative
integers such that every prime number $4p-1$ that divides $k$
occurs even time in the prime factorization of $k$.
\end{conjecture}

Conjecture \ref{conj1} could be reformulated in the following way:
\medskip

\textit{There is an integer matrix $P=[p_{i j}]$ $1\leq\, i, j\leq
4 n+2$ such that $$\sum_{j}^{4n+2} p_{i j}^2 =k$$ for every $i=1,
\dots, 4 n+2$ and $$\sum_{t}^{4n+2} p_{i t} p_{j t}=0$$ for every
$i\neq j$ iff $k$ could be written as the sum of two squares?}
\medskip

We could think about the columns of $P$ as the vectors in
$\mathbb{R}^{4n+2}$. Let observe that matrix $P$ would satisfy the
equality case in the famous Hadamard's Inequality (see
\cite{Bullen}, p. 108). This means that if we look at the columns
of $P$ as generators of the lattice (which is the sublattice of
$\mathbb{Z}^{4n+2}$) the integer $k$ is its discriminant. Our
question is what are the values of discriminants of orthogonal
integer lattices in $\mathbb{R}^{4n+2}$ with equal lengths of
generators. The matrices that satisfy equality case of Hadamard's
Inequality are frequently seen in mathematics. Those with entries
$-1$ and $1$ are called Hadamard's matrices (see \cite{Agaian}).
There are no $(4n+2)\times (4n+2)$ Hadamard's matrices by the
result of Paley from 1933. We think that Conjecture \ref{conj1} is
highly connected with studying the orthogonal lattices and their
discriminants.

\section{Some observations about maps between quasitoric $4$-manifolds}

In the previous sections, we saw several examples of the sets $D
(M, N)$ when $M$ and $N$ are quasitoric $4$-manifolds. We are not
able to determine this set in general for quasitoric $4$-manifolds
but due to Theorem \ref{bar1} and Corollary \ref{bara2} and
special cases from Section 4 we could determine it for various
manifolds, and in most of cases say something about them.

We could think about problem in this manner. We decompose $M$ and
$N$ as the connected sums of $\mathbb{C} P^2$,
$\overline{\mathbb{C} P^2}$ and $S^2\times S^2$. From the system
of Diophantine equations we could notice some general restriction
on degree $k$ map, for example that should be positive or negative
or even or the sum of certain number of squares or no restriction.
Then we are working backward. We studied maps to $\mathbb{C} P^2$,
$\overline{\mathbb{C} P^2}$, $S^2\times S^2$, $\mathbb{C}
P^2\sharp\overline{\mathbb{C} P^2}$, $\mathbb{C}
P^2\sharp\mathbb{C} P^2$ etc. We could combine this results by
repeatedly application of Theorem \ref{bar1} and hope to produce a
degree $k$ map $f: M\rightarrow N$. This is not always possible to
do, but this is the algorithm for generating new examples.

\begin{example} Let $l$, $m$ and $n$ be positive integers such
that $l\geq m+n$, then there is a degree $k$ map $$f: (S^2\times
S^2)^{\sharp l} \rightarrow ({\mathbb{C} P^2})^{\sharp n}\sharp
(\overline{\mathbb{C} P^2})^{\sharp n} \sharp (S^2\times
S^2)^{\sharp m}$$ if and only if $k$ is even number.
\end{example}

\begin{example} Let $m$ and $n$ be positive integers such
that $m\geq n$, then there is a degree $k$ map $$f: ({\mathbb{C}
P^2})^{\sharp 2 m+1} \rightarrow ({\mathbb{C} P^2})^{\sharp
n}\sharp(S^2\times S^2)^{\sharp n}$$ iff $k$ is the square of an
integer.
\end{example}

\begin{example} Let $l$, $m$, $n$ and $p$ be positive integers such
that $p\geq n$ and $l\geq m$, then for every integer $k$ there is
a degree $k$ map
$$f: ({\mathbb{C} P^2})^{\sharp p}\sharp(\overline{\mathbb{C}
P^2})^{\sharp p} \sharp (S^2\times S^2)^{\sharp l} \rightarrow
({\mathbb{C} P^2})^{\sharp n} \sharp(\overline{\mathbb{C}
P^2})^{\sharp n} \sharp (S^2\times S^2)^{\sharp m}$$.
\end{example}

Having this idea in mind we get the following theorem:

\begin{theorem}\label{bar11} Let $M$ be a given quasitoric $4$-manifold. Then
there exist integers $a_0$, $b_0$ and $c_0$ such that for every
integers $a$, $b$ and $c$ such that $a\geq a_0$, $b\geq b_0$ and
$c\geq c_0$, and $$D (({\mathbb{C} P^2})^{\sharp
a}\sharp(\overline{\mathbb{C} P^2})^{\sharp b} \sharp (S^2\times
S^2)^{\sharp c}, M)= \mathbb{Z}.$$
\end{theorem}

This theorem states that there are infinitely many manifolds that
could be mapped to $M$ with any degree.

\section{Concluding remarks}

In the previous sections we made the numerous calculations. We
used the number theory for determination of the mapping degree
between special class of manifolds. According to the result of
Duan and Wang all these observations are relevant at least as the
necessary condition in more general case of $(n-1)$-connected $2n$
manifolds. One direction for further research is to connect our
result with higher dimensional manifolds that have the same
intersection forms. Duan and Wang also gave sufficient condition
that one need to check for the existence of degree $k$ map. Our
results showed that even in the simplest case of $4$-manifolds, it
is hopeless to directly check conditions.

If we stay on $4$-manifolds, the famous result \cite{Freed} of
Freedman gives classifications of $4$-manifolds in the terms of
their intersection forms. We have already noticed the connection
of the problem with the problems about lattices and quadratic
forms. So it is naturally to study problem from this point of
view. We expect that progress in any of this area may lead to
progress in our problem and vice versa.

The partially results we got in the case of quasitoric
$4$-manifolds give hope to fully answer the question for this
class of manifolds. We think that other more sophisticated
technics of the number theory could at least enlarge our knowledge
about the mapping degrees.

\begin{center}\textmd{Acknowledgements }
\end{center}

\medskip The author is grateful to Rade \v{Z}ivaljevi\'{c} and
Vladimir Gruji\'{c} for fruitful discussion and comments.

\bigskip

{\small \DJ{}OR\DJ{}E BARALI\'{C}, Mathematical Institute SASA,
Kneza Mihaila 36, p.p.\ 367, 11001 Belgrade, Serbia

E-mail address: djbaralic@mi.sanu.ac.rs

\end{document}